\documentclass[a4paper]{article}
\usepackage{cite}
\usepackage{amsmath,amssymb,amsfonts}
\usepackage{algorithm,algorithmic}
\usepackage[dvipdfmx]{graphicx}
\usepackage{url}
\usepackage{textcomp}
\usepackage{xcolor}
\def\BibTeX{{\rm B\kern-.05em{\sc i\kern-.025em b}\kern-.08em
  T\kern-.1667em\lower.7ex\hbox{E}\kern-.125emX}}
\begin{document}

%---------------------------------------%
% Title
%---------------------------------------%
\title{Acceleration of multiple precision solver for ill-conditioned algebraic equations with lower precision eigensolver}
%\titlerunning{Acceleration of MPF solver of algebraic equation}

%---------------------------------------%
% Author
%---------------------------------------%
\author{Tomonori Kouya\thanks{Shizuoka Institute of Science and Technology}\\\url{https://na-inet.jp/}}%\inst{1}\orcidID{0000-0003-0178-5519}}
%
%\authorrunning{T.Kouya}
%\institute{Shizuoka Institute of Science and Technology,\\2200-2 Toyosawa, Fukuroi 437-8555, Japan \\
%\email{kouya.tomonori@sist.ac.jp}}
\date{March 28, 2022}

\maketitle

%---------------------------------------%
% Abstract
%---------------------------------------%
\begin{abstract}
There are some types of ill-conditioned algebraic equations that have difficulty in obtaining accurate roots and coefficients that must be expressed with a multiple precision floating-point number. When all their roots are simple, the problem solved via eigensolver (eigenvalue method) is well-conditioned if the corresponding companion matrix has its small condition number. However, directly solving them with Newton or simultaneous iteration methods (direct iterative method for short) should be considered as ill-conditioned because of increasing density of its root distribution. Although a greater number of mantissa of floating-point arithmetic is necessary in the direct iterative method than eigenvalue method, the total computational costs cannot obviously be determined. In this study, we target Wilkinson's example and Chebyshev quadrature problem as examples of ill-conditioned algebraic equations, and demonstrate some concrete numerical results to prove that the direct iterative method can perform better than standard eigensolver.
%\keywords{multiple precision floating-point arithmetic \and algebraic equation \and eigensolver \and parallelization}
\end{abstract}

%--------------------------------------%
% Introduction
%--------------------------------------%
\section{Introduction}

Currently, multiple-precision floating-point (MPF for short) arithmetic are executed using reliable and highly-performed de-facto standard libraries such as QD\cite{qd} and MPFR\cite{mpfr}. These have been developed since the end of the 20th century, and have been in use for over two decades. More convenient and efficient multiple precision numerical computation libraries such as MPLAPACK\cite{mplapack} and our BNCpack\cite{bnc}, are constructed on that of MPF libraries. Therefore, we can easily deal with various types of ill-conditioned problems in normal hardware and software environment for consumers.

We have published a paper\cite{kouya_gauss_en} about deriving highly accurate abscissas of Gauss-type integration formulas using both "eigenvalue method" via symmetric tridiagonal matrix originated by Golub and Welsch\cite{golub_welsch} and "direct iterative method" to directly calculate zeros of orthogonal polynomial using Newton iteration\cite{yamashita_gauss_legendre_en}\cite{yamashita_gauss_laguerre_en}\cite{yamashita_gauss_hermite_en}. Our previous task with direct iterative method was to calculate the approximation of abscissas with user-required accuracy by combining binary64 eigensolver of LAPACK\cite{lapack} and MPFR Newton iteration. 

The performance of the eigenvalue method is slow but robust for effect by round-off error; hence, we can obtain good approximations even with binary64 arithmetic. However, the direct iterative method is well-performed with multiple precision arithmetic, but good initial guesses must be employed to guarantee its convergence. In our previous paper, we concluded that, to shorten total computational time, it is best to combine binary64 LAPACK eigensolver to derive initial guesses, and multiple precision Newton iteration to exploit a higher accuracy. From the current point of view, our proposed combination employs the ``mixed precision" method for solving algebraic equations.

From these experiences, we have tried to combine multiple-precision eigensolver of MPLAPACK, and second and third order simultaneously iterative methods of BNCpack, such as Durand-Kerner (DK for short) method, to solve two types of ill-conditioned algebraic equations with real and complex roots. One is well-known Wilkinson's problem\cite{wilkinson_roundoff}, and another is calculation of abscissas of Chebyshev quadrature formula (Chebyshev quadrature problem for short), which has been experimented by Harumi Ono\cite{ono_dka1_en}\cite{ono_dka2_en}. The algebraic equation of Chebyshev quadrature problem has already been studied in several published papers, but those papers are not known except in Japan because almost these papers are written only in Japanese. The only English non-referred paper in English can be found in Kokyoroku series of RIMS in Kyoto University\cite{moriguchi_iri_ono_en}. Here, we briefly explain the previous results on the Chebyshev quadrature problem.

The algebraic equations derived from the Chebyshev quadrature formula have difficulties not only in calculating the accurate coefficients, but also in solving the equations with multiple precision floating-point arithmetic. Harumi Ono published papers on how to solve the 1024-th degree polynomial in 1979, and 2048-th degree in 1981 using DK method with her original multiple precision arithmetic on Cray. In addition, Masumoto et.al. reported the numerical property of 20480-th degree coefficients when deriving them with MPF arithmetic\cite{ono_fujino_en}, but did not achieve the obtained roots of those. We will describe the mathematical expression of the Chebyshev quadrature problem in Section \ref{section:examples_illcondition}.

In this study, we will relay the results obtained by combining DD Rgeev of MPLAPACK as a generator of initial guesses and the MPFR Durand-Kerner method from the point of view for mixed precision approach of acceleration.

The following two computational environments, EPYC and Xeon, are used in the rest of this study. MPLAPACK and our BNCpack, including QD and MPFR/GNU MP, are natively compiled with Intel Compiler.
\begin{description}
	\item[EPYC] AMD EPYC 7402P 24 cores, Ubuntu 18.04.6 LTS, GCC 7.5.0, Intel Compiler version 2021.4.0, MPLAPACK 1.0.1, BNCpack 0.8, MPFR 4.1.0
	\item[Xeon] Intel Xeon W-2295 3.0GHz 18 cores, Ubuntu 20.04.3 LTS, GCC 9.3.0, Intel Compiler version 2021.5.0, MPLAPACK 1.0.1, BNCpack 0.8, MPFR 4.1.0
\end{description}

OpenMP is applied for parallelization of DK method with the following compile option as \verb|icpc -O3 -qopenmp|.

%--------------------------------------%
% 
%--------------------------------------%
\section{Belief introduction of current multiple precision arithmetic and performance of current MPFR GEMM}

The need for computational processing, mainly using floating-point arithmetic, will increase and not decrease in the foreseeable future, not only for scientific computing, but also for deep learning and other applications. The performance of hardware such as CPUs and GPUs is mainly improved in the following two ways:
\begin{itemize}
  \item Parallel processing capabilities such as SIMD instructions and an increase in the number of cores in CPUs and GPUs,
  \item Increasing or decreasing the length of mantissa and exponent in floating-point numbers according to the necessity of user's computation, such as the introduction of half-precision and single-precision floating-point operations.
\end{itemize}

In similar ways, the whole process of computation is now accelerated by enhancing parallelization techniques. While taking advantage of these hardware performance improvement and employing parallelization techniques such as SIMD instruction, OpenMP, and MPI, mixed precision techniques combining IEEE binary32, binary64, and half precision arithmetic, which are standard in hardware, are also used. Regarding MPF operations using software libraries, which are adopted in adverse conditions where binary64 arithmetic lacks computational accuracy, the heavy MPF processing requires the active use of similar software performance improvement techniques.

Currently, the mainstream of MPF arithmetic falls into two types of implementations: a multi-component method that combines multiple binary32 and binary64, using error-free transformation (EFT) techniques to extend the mantissa length, and an integer-based many-digit method. The QD library by Bailey et al. is well-known for those based on the multi-component method, and the MPFR library using the arbitrary-length natural number kernel (MPN) of GNU MP (GMP) has a significant number of users because of its superior speed and reliability. MPLAPACK by Maho Nakata is a multiple precision linear computation library based on the C++ converted from the original Fortran code of LAPACK/BLAS, which was developed from these two multiple precision floating-point arithmetic libraries, and the latest version (Version 1.0.1) as of February, 2022, provides the parallelized BLAS code with OpenMP. In addition, the main driver and calculation routines of LAPACK are available in multiple precision. However, we will wait to gain all driver routines benefit from this high-performance MPBLAS. In addition, SIMD instructions, the use of CUDA, and the introduction of the Ozaki scheme have not officially improved performance as faster multiple precision ATLAS, OpenBLAS, and Intel Math Kernel.
	
\begin{figure}[htbp]
	\begin{center}
		\includegraphics[width=.45\textwidth]{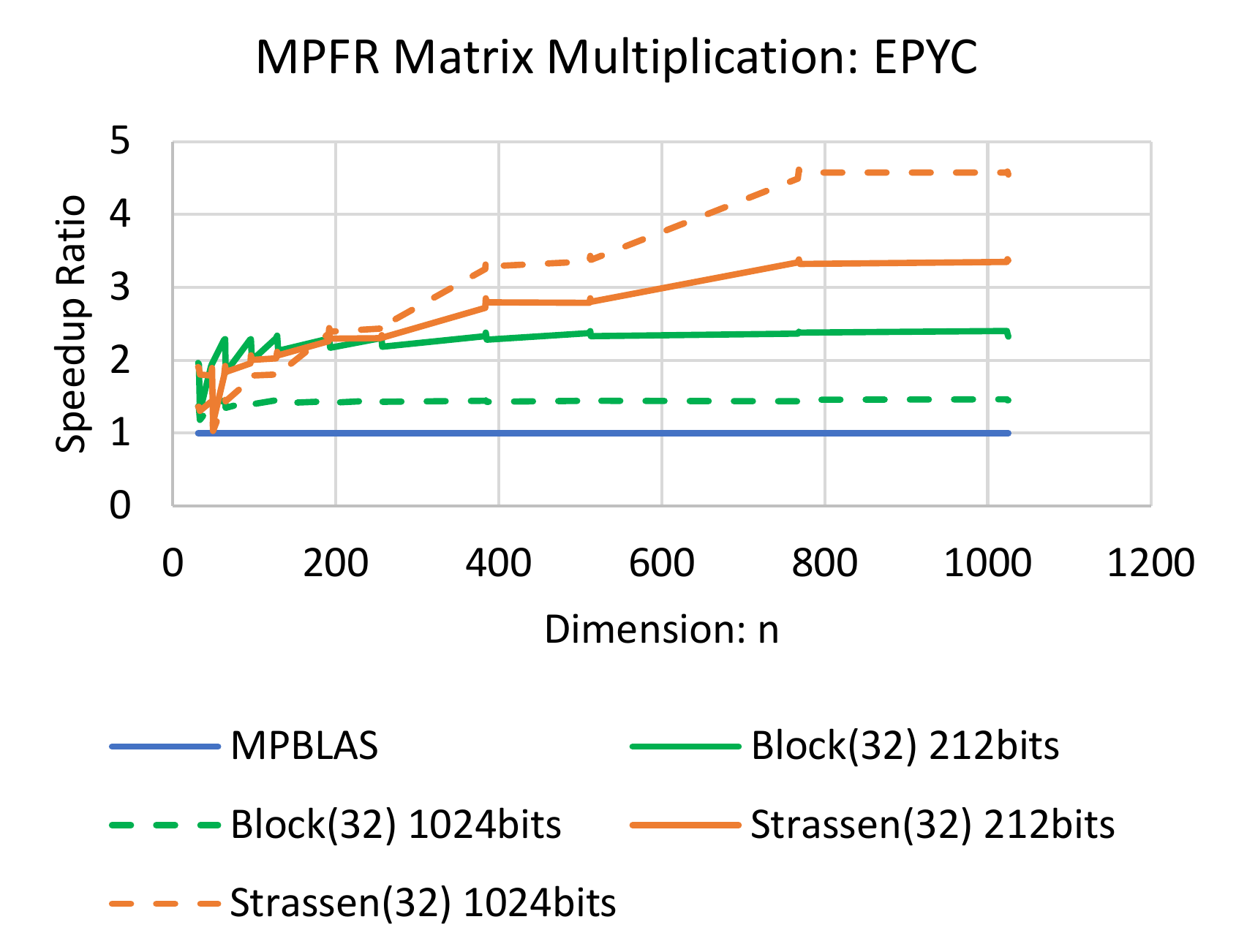}
		\includegraphics[width=.45\textwidth]{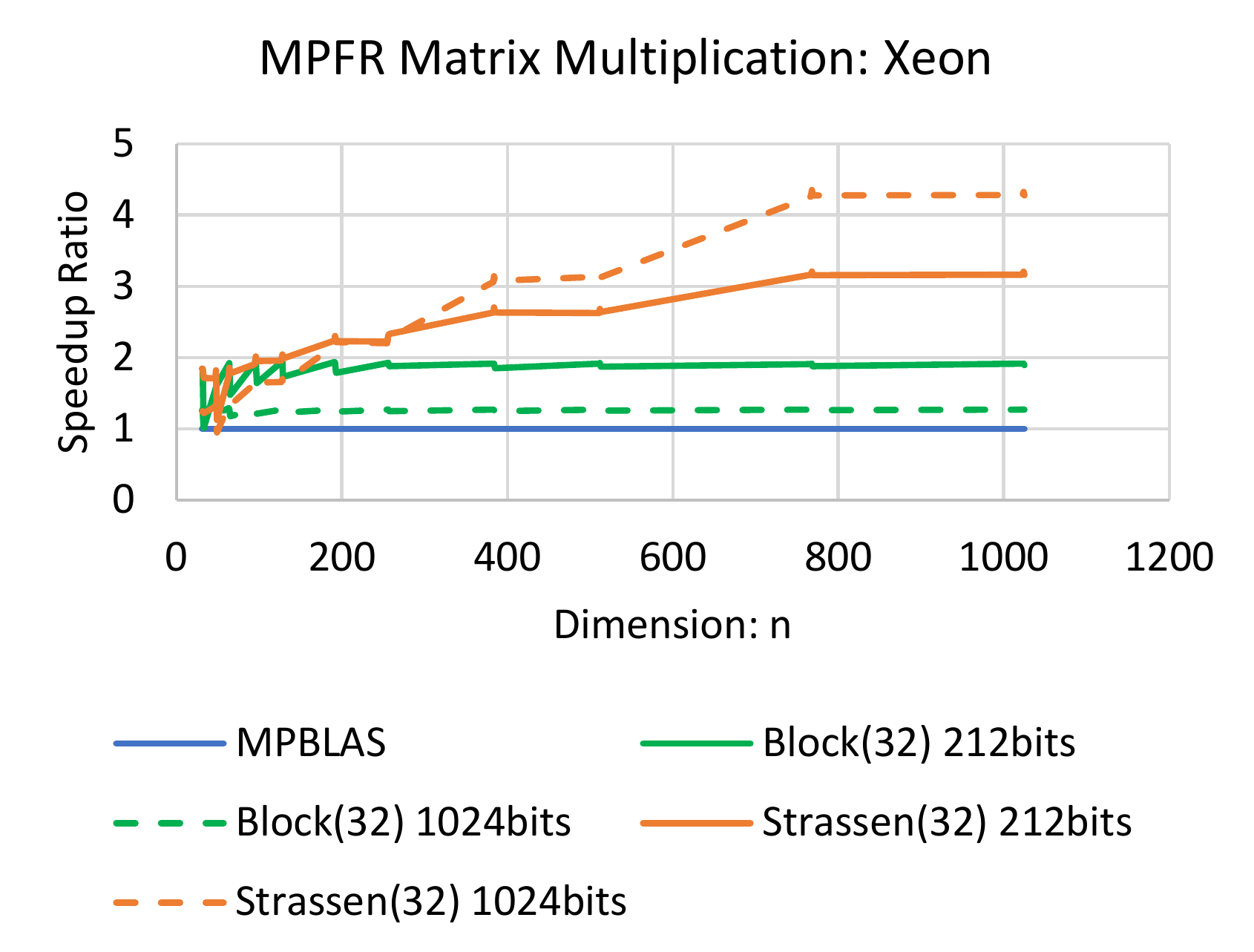}
		\caption{Speedup ratio of MPFR matrix multiplication against MBLAS(Rgeem)}\label{fig:matmulbench_mpfr}
	\end{center}
\end{figure}

The current driver routine Rgeev of MPLAPACK used in this paper is not parallelized. In addition, the arbitrary-precision routines have been determined to be slower than our native C implementation of the basic linear subprogram, probably because of the use of MPREAL, a wrapper C++ class library for MPFR adopted in the arbitrary-precision calculations. In fact, the computational time of the MPFR block matrix multiplication supported by our BNCmatmul library is illustrated in \figurename\ref{fig:matmulbench_mpfr} compared to the time of the MPFR Rgeem of MPBLAS. The results of Strassen matrix multiplication are also included for comparison.

Evidently, both 212 bits (the same bits of QD precision) and 1024 bits are 2.0 – 2.4 times (212 bits) and 1.4 – 1.5 times (1024 bits) larger than MPBLAS (Rgemm) by block matrix multiplication (matmul\_mpfmatrix\_block) in both EPYC and Xeon environments. Thus, the smaller the number of MPF mantissa, the faster the C native MPFR direct call (BNCmatmul) is than MPLAPACK/MPBLAS (MPREAL).

The eigenvalue driver routine (Rgeev) for real matrices in MPLAPACK are adopted as eigensolvers, while the simultaneous iterative method for MPFR algebraic equations, prepared as a competitor, employs a direct call to MPFR for its implementation. Therefore, it is expected to perform better than the implementation using MPREAL.

%--------------------------------------%
% 
%--------------------------------------%
\section{Mixed precision approach for solving algebraic equations}

We target the following algebraic equations
\begin{equation}
  p_n(x) = 0, \label{eqn:algebraic_eq}
\end{equation}
where the following real coefficient polynomial $p_n(x)$ is adopted as the left term,
\begin{equation}
  p_n(x) = \sum^n_{i=0} a_i x^i\ (a_i\in\mathbb{R}, a_n\not=0). \label{eqn:pn}
\end{equation}

It is well-known that the $n$-th degree algebraic equation (\ref{eqn:algebraic_eq}) definitely has $n$ roots $\alpha_i \in \mathbb{C}$, $i = 1, 2, ..., n$ at most. In this time, suppose that we have no prerequisite for roots of the equation (\ref{eqn:algebraic_eq}).

For convenience, the corresponding monic polynomial $q_n(x)$ derived from $p_n(x)$
\begin{equation}
  q_n(x) = x^n + \sum^{n-1}_{i=0} c_i x^i\ (c_i = a_i/a_n), \label{eqn:qn}
\end{equation}
is also prepared.

%--------------------------------------%
% 
%--------------------------------------%
\subsection{Eigenvalue method for roots of algebraic equation}

As explained in several number of textbooks of linear algebra, it is well-known that our targeted algebraic equation (\ref{eqn:algebraic_eq}) can be expressed by the following companion matrix $C_n$ with the same eigenvalues as the roots of $q_n(x)$. %The $C_n$ is expressed as 
\begin{equation}
	C_n = \left[\begin{array}{ccccc}
		0   & 1   & 0   & \cdots  & 0 \\
		\vdots & \ddots & \ddots & \ddots  & \vdots \\
		0   & \cdots & 0   & 1    & 0 \\
		0   & \cdots & \cdots & 0    & 1 \\
		-c_0  & -c_1  & \cdots & -c_{n-2} & -c_{n-1}
	\end{array}\right]. \label{eqn:companion}
\end{equation}

As the polynomial (\ref{eqn:pn}) has only real coefficients, we can exploit xGEEV driver routines in LAPACK\cite{lapack} as an eigensolver to obtain all the eigenvalues of $C_n$. If $C_n$ has a negligible condition number, is diagonalizable, and does not have multiple eigenvalues, we do not require additional digits of mantissa of MPF. We call this approach the “eigenvalue method” for solving the algebraic equation (\ref{eqn:algebraic_eq}).

%--------------------------------------%
% 
%--------------------------------------%
\subsection{Simultaneous iteration method for directly solving algebraic equation} 

The standard direct iterative method for obtaining solutions to algebraic equations of degree five or higher (\ref{eqn:algebraic_eq}) is the Newton method and its related simultaneous methods. In recent years, several higher order direct iterative methods have been proposed, and Petkovic has compactly summarized the results up to 2012\cite{petkovic_book}. However, as described below, the higher the order of the direct iterative method, the more computational complexity per iteration increase, and global convergence is not guaranteed. Therefore, we employ only the second- and third- order methods\cite{aberth_3rd}. The comparison of the initial guess setting methods will be discussed later.

For both second- and third- order DK methods, by expressing the approximation at $k$ times iteration as 
\[ \mathbf{z}_k = [z^{(k)}_1\ z^{(k)}_2\ ...\ z^{(k)}_n]^T\in \mathbb{C}^n, \]
their iteration formulas are described with monic polynomial (\ref{eqn:qn}) as follows:
\begin{description}
	\item[Second order DK Method]
\begin{equation}
  z^{(k+1)}_i := z^{(k)}_i - \frac{q_n(z^{(k)}_i)}{\prod^n_{j=1, j\not=i} (z^{(k)}_i - z^{(k)}_j)}\label{eqn:dk2}
\end{equation}
	\item[Third order DK method]
\begin{equation}
  z^{(k+1)}_i := z^{(k)}_i - \frac{\frac{p_n(z^{(k)}_i)}{p_n'(z^{(k)}_i)}}{1-\frac{p_n(z^{(k)}_i)}{p_n'(z^{(k)}_i)}\displaystyle\sum^n_{\substack{j=1\\j\not=i}} (z^{(k)}_i - z^{(k)}_j)^{-1}}\label{eqn:dk3}
\end{equation}
\end{description}

To compare with our initial guess approach, we adopt Aberth's initial guesses as follows:
\begin{equation}
	z^{(0)}_i := -\frac{c_{n-1}}{n} + r\exp\left\{\left(\frac{2(i-1)\pi}{n}+\frac{3}{2n}\right)\mathrm{i}\right\}. \label{eqn:dka_init}
\end{equation}

We adopt $r$ in Aberth's initial guess (\ref{eqn:dka_init}) as follows:
\[ r := \max_{0\le i \le (n-1)} |n_{\rm nz} c_i|^{1/(n - i)}, \]
where $n_{\rm nz} \le n$ is the number of non-zero coefficients in (\ref{eqn:pn}).

%--------------------------------------%
% 
%--------------------------------------%
\section{Examples of ill-conditioned algebraic equations}\label{section:examples_illcondition}

Unlike the eigenvalue problem, when solving the equation, the higher the density of roots, the larger the error in the values of the polynomial $p_n(x)$ and $q_n(x)$. Therefore, if we employ a method such as the Danilevsky method \cite{fadeev_en}, where the eigenvalue problem of a matrix is transformed into a companion matrix and coefficients of the eigenequation are obtained directly, despite the good conditions for an eigenvalue problem, it becomes a bad problem to solve the algebraic equation. This is why this method via eigenequation is not recommended currently because fast eigenvalue solving methods such as the QR method, which uses a shift of the origin in the Householder transformation to obtain stable and accurate eigenvalues exist.

However, regarding bad algebraic equations that require longer MPF numbers than binary64 for the coefficients, MPF operations are essential in the process of solving, and there is a possibility that the solution of algebraic equations with less computational complexity can be performed faster than matrix eigenvalue routines. In our previous study on the Gauss-type integral formulas for the quantile calculations \cite{kouya_gauss_en}, there were several cases where the computational time was reduced using the Newton method with the approximate eigenvalues obtained by binary64 as initial guess, rather than computing all eigenvalues using the MPF eigensolver for real symmetric matrices. As an extension of this result, it is expected that similar mixed-precision techniques will be effective in reducing computational time for solving arbitrary real coefficient algebraic equations using a real asymmetric companion matrix to derive useful initial guesses.

Here we consider two examples: Wilkinson's example (with all real roots) and an algebraic equation whose solution is the quantile of the Chebyshev quadrature formula (with almost all complex roots), to prove that such examples exist in the benchmark test.

%--------------------------------------%
% 
%--------------------------------------%
\subsection*{Wilkinson's example}

This well-known example is written and explained in detail in Wilkinson's book\cite{wilkinson_roundoff}.

\begin{figure}[htb]
	\begin{center}
	\includegraphics[width=.6\textwidth]{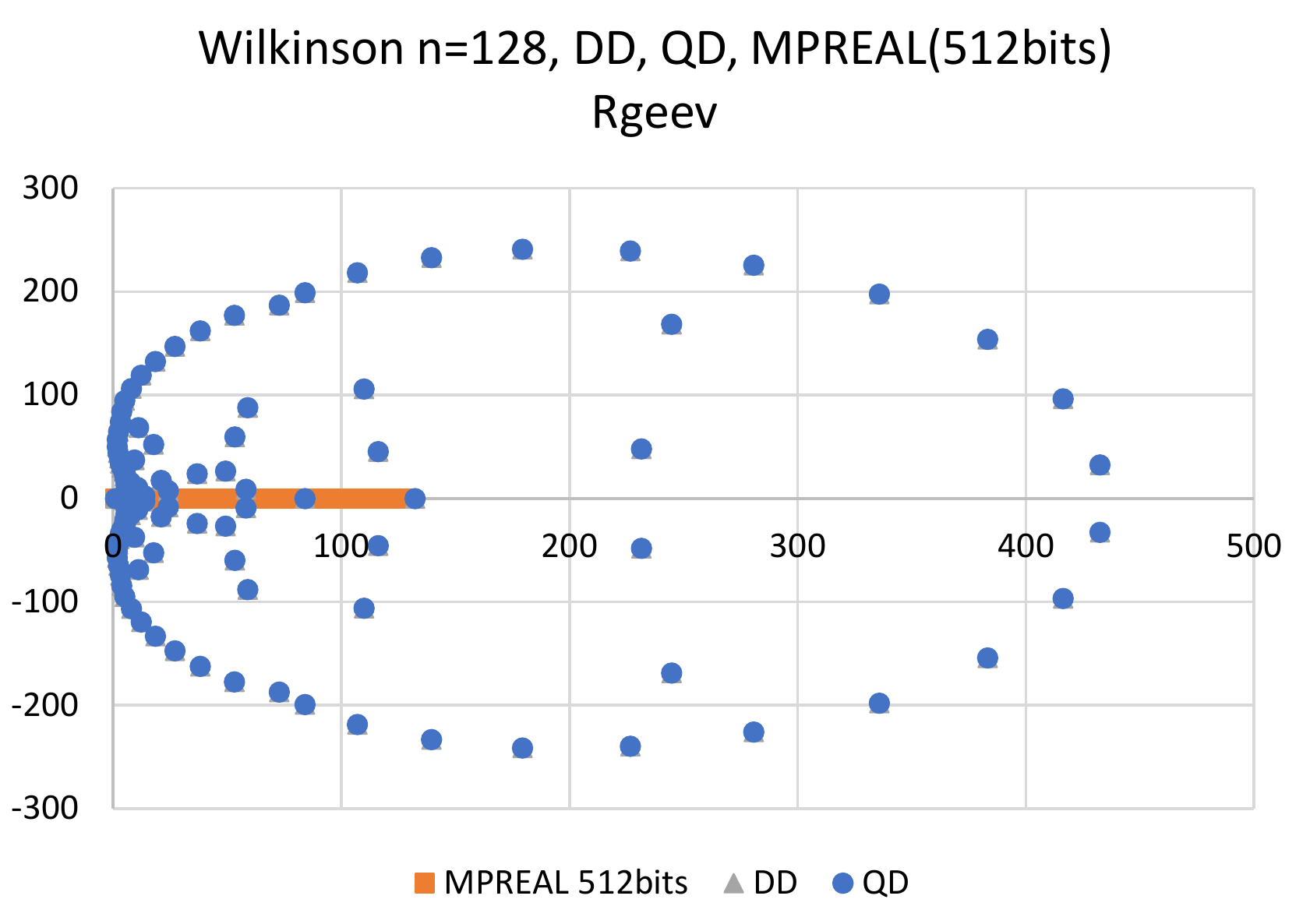}
	\caption{DD, QD, MPREAL(512 bits) Rgeev(MPLAPACK): Wilkinson's example $n=128$}\label{fig:wilkinson128deg}
	\end{center}
\end{figure}

The polynomial (\ref{eqn:pn}) is provided as $p_n(x) = \prod^n_{i=1} (x - i)$, with $\alpha_i = i$ as its roots. The absolute value of coefficients is glowing up on $O(n!)$. For $n=20$, the coefficients are
	\begin{eqnarray*}
		a_0 &=& 2432902008176640000\\
		a_1 &=& -8752948036761600000\\
		%a_2 &=& 13803759753640704000\\
		&\vdots& \\
		%a_{18} &=& 20615\\
		a_{19} &=& -210\\
		a_{20} &=& 1.
	\end{eqnarray*}
In the case of $n=128$ adopted in our experiments, $a_0 = 3.8562\cdots\times 10^{215}$ means that coefficients at over $n=128$ may not be treated in binary64, DD (Double-double, 106 bits), and QD (Quad-double, 212 bits) arithmetic. \figurename\ref{fig:wilkinson128deg} illustrates the approximation of the roots with DD, QD, and MPREAL 512 bits Rgeev. It is obvious that DD and QD Rgeev cannot provide an accurate approximation of the roots.

However, as indicated later, DD approximations are useful as initial guesses to accelerate MPF direct iterative methods.

%--------------------------------------%
% 
%--------------------------------------%
\subsection*{Chebyshev quadrature problem}

As previously described, this example is solely popular in Japan. In this subsection, we mathematically explain the Chebyshev quadrature problem.

The Chebyschev quadrature formula is not categorized in the Gaussian quadrature formula. When $[-1, 1]$ is provided as the integration interval, we can obtain the following discrete quadrature formula as follows:
\begin{equation}
	\int^{1}_{-1} w(x) f(x) dx \approx \sum^n_{k=1} w_k(x_k) f(x_i). \label{eqn:chebyint}
\end{equation}
In the above formula (\ref{eqn:chebyint}), we suppose that $w(x)=1$ and $w_k(x_k) = 2/n$ is fixed. In this case, the “Chebyshev quadrature formula” is determined by choosing the appropriate abscissas.

The coefficients of polynomial (\ref{eqn:pn}) with roots as their abscissas are derived, starting with $a_n=1$, as follows:
\begin{equation}
\begin{cases}
	a_{n-(2k-1)} := 0 & \\
	a_{n - 2k} := -\sum^k_{2j+1} a_{n-2(k-j)} & \\
\end{cases}, \label{eqn:chebcoef}
\end{equation}
where $k = 1, 2, ..., \lfloor n/2 \rfloor$. The odd numbered terms are zero, and the even numbered terms are derived from (\ref{eqn:chebcoef}).

The algebraic equation (\ref{eqn:algebraic_eq}) with these above coefficients has real roots only in the case of $n = 1, 2, ..., 7$, $9$, and almost all conjugate complex pairs of roots in the case of $n = 8$, or $ n \geq 10$. According to Moriguchi and Iri's prediction\cite{moriguchi_iri_ono_en}, regarding $n\rightarrow \infty$, all roots are nearly distributed on the curve and expressed as the following functions of $z$,
\begin{equation}
	|(z + 1)^{(z+1)/2} (z-1)^{-(z-1)/2}| = 2. \label{eqn:cheby_roots}
\end{equation}
The prediction has been precisely confirmed in numerical experiments through Ono's works \cite{ono_dka1_en}\cite{ono_dka2_en}.

\begin{figure}[htbp]
	\begin{center}
	\includegraphics[width=.6\textwidth]{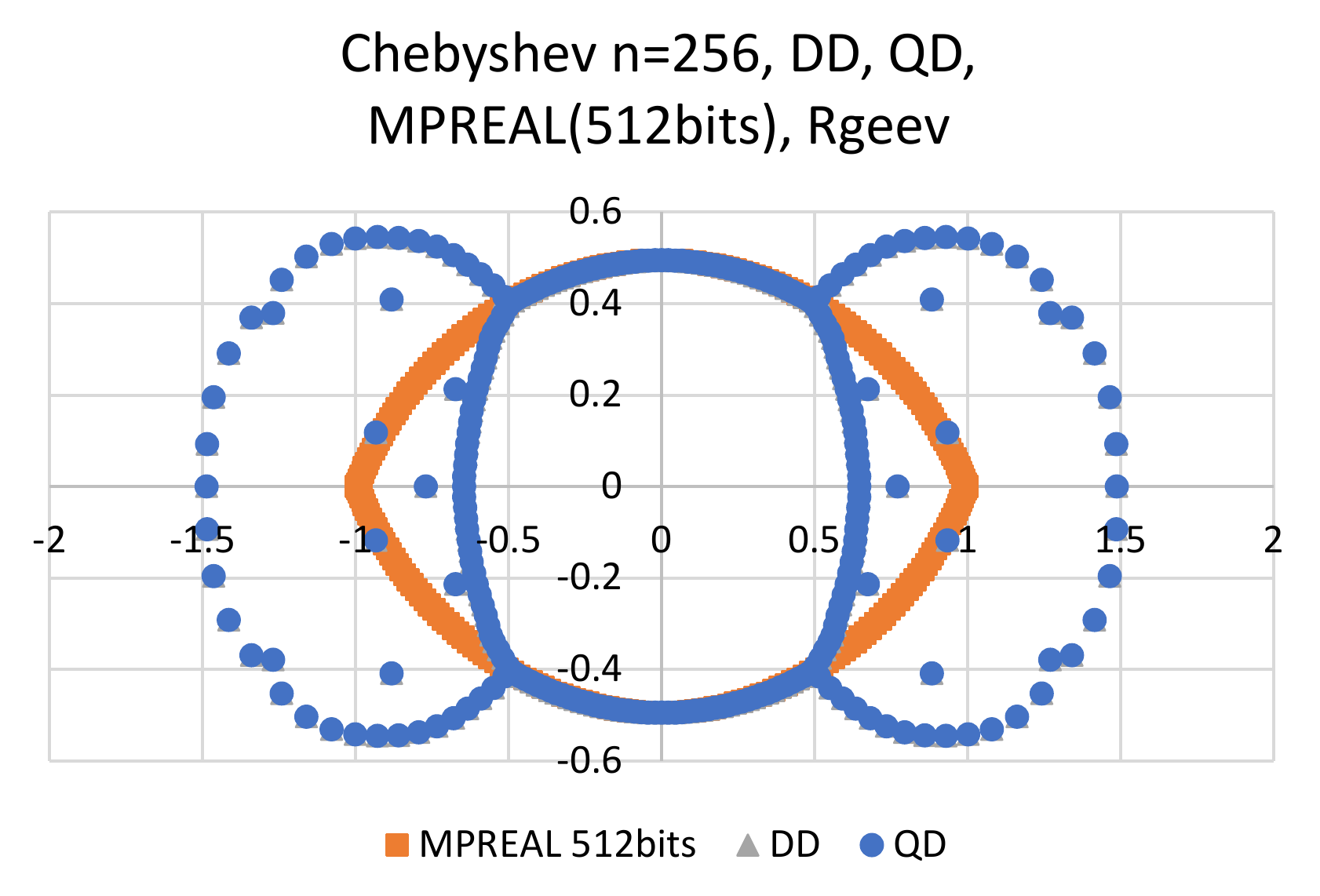}
	\caption{The abscissas of Chebyshev quadrature problem ($n = 256$) obtained by DD, QD, MPREAL(512 bits) Rgeev of MPLAPACK}\label{fig:ono256deg}
	\end{center}
\end{figure}

In addition, we should handle the process of calculating the coefficients using the formula ($\ref{eqn:chebcoef}$). With an increase in the order of the polynomial, the digits decrease significantly, and longer MPF operations are required to obtain accurate coefficients. According to studies of Masumoto et.al., we need over 215 decimal digits to guarantee the over six decimal significant digits of $a_0$ in case of $n=1024$. As illustrated in \figurename\ref{fig:ono256deg}, DD and QD precision arithmetic are not sufficient to obtain accurate roots even in the case of $n=256$.

%--------------------------------------%
% 
%--------------------------------------%
\section{Benchmark tests}

Our mixed precision approach for targeted algebraic equation (\ref{eqn:algebraic_eq}) is set up as follows:
\begin{enumerate}
	\item Calculate all coefficients $a_i$, $i=0, 1, ..., n$ of the algebraic equation with $L$ bits MPFR arithmetic.
	\item Transfer MPFR $a_i$ to the DD ones and then calculate only eigenvalues $\lambda_i^{(\mathrm{DD})}\in \mathbb{C}$ of the corresponding companion matrix $C$ using DD Rgeev of MPLAPACK.
	\item Set the above DD eigenvalues $\lambda_i^{(\mathrm{DD})}$ as initial guesses $\mathbf{z}_0\in\mathbb{C}^n$ for second- or third- order DK methods. Converged values $\mathbf{z}_{\rm end} \in$ $\mathbb{C}^n$ are adopted as final approximations of roots $\alpha_i$.
\end{enumerate}

Although the value of $|p_n(x)|$ is normally used for the stopping rule of iteration, we use the difference of $z_i^{(k)}$ and $z_i^{(k+1)}$, relative tolerance $\varepsilon_{\rm rel}$, and absolute tolerance $\varepsilon_{\rm abs}$ as follows:
\begin{equation}
	|z_i^{(k+1)} - z_i^{(k)}| \leq \varepsilon_{\rm rel} |z_i^{(k)}| + \varepsilon_{\rm abs}. \label{eqn:stopping_rule}
\end{equation}
We recognize that the $(k+1)$-th iterated approximation is converged when the condition of (\ref{eqn:stopping_rule}) is satisfied. The relative errors as illustrated in latter figures are calculated by comparing with the eigenvalues using MPREAL 2048 bits Rgeev.

%--------------------------------------%
% 
%--------------------------------------%
\subsection{Wilkinson's example}

First, we relay the results for Wilkinson's example ($n = 128$) on EPYC and Xeon.

\begin{figure}[htbp]
	\begin{center}
		\includegraphics[width=.55\textwidth]{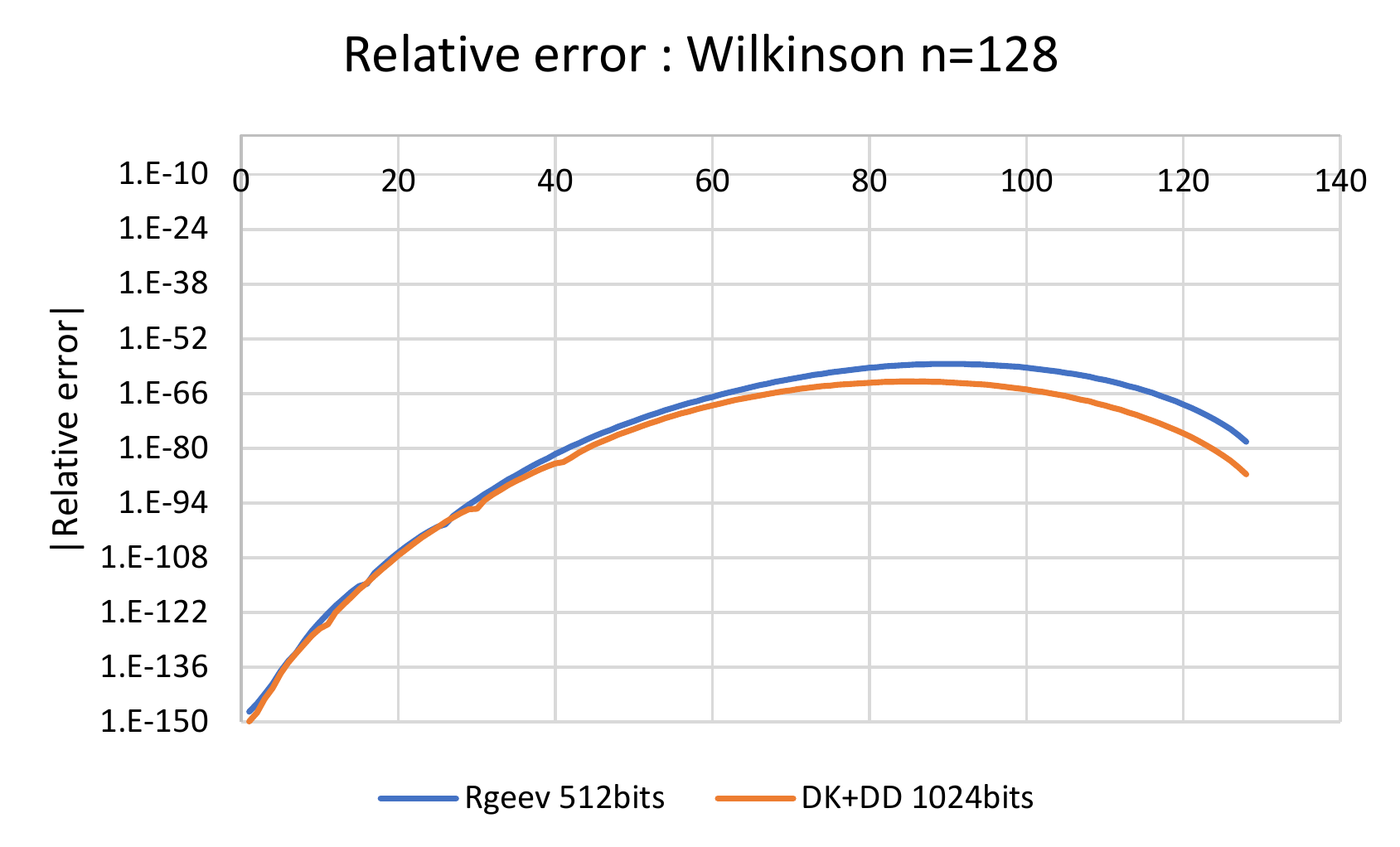}
		\caption{Relative errors of Rgeev(512 bits) and DK-DD(1024 bits) for Wilkinson's example: $n=128$.}\label{fig:relerr_wilkinson_128deg_rgeev512_dka1024}
	\end{center}
\end{figure}

The 512 bits MPREAL Rgeev can obtain an accurate approximation of roots. In contrast, to obtain an approximate solution with sufficient accuracy using the DK method, the 1024 bits are necessary, and $\varepsilon_{\rm rel} := 7.5\times 10^{-145}$ and $\varepsilon_{\rm abs} := 1.0\times 10^{-300}$ are adopted to determine convergence. Accordingly, we can confirm that the approximations of both ways can reach the same level of relative errors as illustrated by \figurename\ref{fig:relerr_wilkinson_128deg_rgeev512_dka1024}. 

The entire results via benchmark tests are presented in \tablename\ \ref{table:wilkinson128}, where “DKA2” (second order DK) and “DKA3” (third order DK) indicate the computational time (in seconds) and number of iterations for each number of threads when the initial guess of Aberth is employed, and “DK2+DD” and “DK3+DD” indicate the computational time and number of iterations when DD Rgeev is adopted to derive the initial guess. Currently, the Rgeev of MPLAPACK is not parallelized; hence, it cannot be accelerated with over two threads. Underlined computational time means being faster than those of MPFR Rgeev.

\begin{table}[htpb]
\begin{center}\small
	\caption{Wilkinson's example:$n=128$, Rgeev 512 bits, DK methods 1024 bits}\label{table:wilkinson128}

	\scalebox{0.8}{\begin{tabular}{|c|c|c|c|c|c|c|c|c|c|c|}
	\multicolumn{1}{c|}{EPYC} & \multicolumn{2}{|c|}{Rgeev (Second)} &\multicolumn{2}{|c|}{DKA2(1024 bits)}	&\multicolumn{2}{|c|}{DKA3(1024 bits)}	&\multicolumn{2}{|c|}{DK2(1024 bits)+DD}	&\multicolumn{2}{|c|}{DK3(1024 bits)+DD}	\\ \hline
		\#Thr.&	DD	&MPFR 512 bits 	& Second	&\#Iter.&Second	&\#Iter.&Second	&\#Iter.&Second	&\#Iter. \\ \hline
	1 	&0.106	&5.1	&79.7			&1374	&71.1	&691	&30.8	&538	&51.5	&512 \\
	2 	&		&		&39.8			&1374	&36.6	&691	&15.4	&535	&26.1	&509 \\ 
	4 	&		&		&20				&1374	&18.7	&691	&8.14	&562	&14		&512 \\
	8 	&		&		&10.1			&1374	&8.96	&691	&\underline{1.43}	&195	&\underline{1.42}	&99 \\
	16	&		&		&5.15			&1374	&\underline{4.55}	&691	&\underline{2.07}	&557	&\underline{3.32}	&511 \\
	24	&		&		&\underline{3.88}	&1374	&\underline{3.74}	&691	&\underline{1.52}	&541	&\underline{2.69}	&511 \\ \hline
	\end{tabular}}
	
	\scalebox{0.8}{\begin{tabular}{|c|c|c|c|c|c|c|c|c|c|c|}
	\multicolumn{1}{c|}{Xeon} & \multicolumn{2}{|c|}{Rgeev (Second)} &\multicolumn{2}{|c|}{DKA2(1024 bits)}	&\multicolumn{2}{|c|}{DKA3(1024 bits)}	&\multicolumn{2}{|c|}{DK2(1024 bits)+DD}	&\multicolumn{2}{|c|}{DK3(1024 bits)+DD}	\\ \hline
	\#Thr.&	DD	&MPFR 512 bits 	& Second	&\#Iter.&Second	&\#Iter.&Second	&\#Iter.&Second	&\#Iter. \\ \hline
	1 & 0.09	& 3.5			& 58.6	&1374	&50.8	&691	&22.3	&538	&37.1	&512\\
	2 & 		& 				& 29.2	&1374	&26.2	&691	&11.3	&535	&18.8	&509\\
	4 & 		& 				& 15.3	&1374	&13.7	&691	&6.18	&562	&9.93	&512\\
	8 & 		& 				& 7.78	&1374	&6.86	&691	&\underline{1.1}	&195	&\underline{0.996}	&99\\
	16 & 		& 				& 4.73	&1374	&4.16	&691	&\underline{1.86}	&557	&\underline{2.97}	&511\\
	18 & 		& 				& 4.88	&1374	&4.37	&691	&\underline{1.9}	&557	&\underline{3.05}	&511 \\ \hline
\end{tabular}}
\end{center}
\end{table}

We can observe the following results from \tablename\ \ref{table:wilkinson128}:

\begin{enumerate}
	\item When Aberth's initial values are adopted, the number of iterations for DK3 is approximately half that of DK2. However, the computational time has not decreased much.
	\item When using the eigenvalues of DD Rgeev as initial values, the number of iterations is reduced compared to when using Aberth's initial values, and the computational time is also reduced. In addition, the number of iteration varies from thread to thread, and the change in $\lambda_i^{(\mathrm{DD})}$ has a significant effect on the iterative process.
\end{enumerate}

\begin{figure}[htbp]
	\begin{center}
		\includegraphics[width=.45\textwidth]{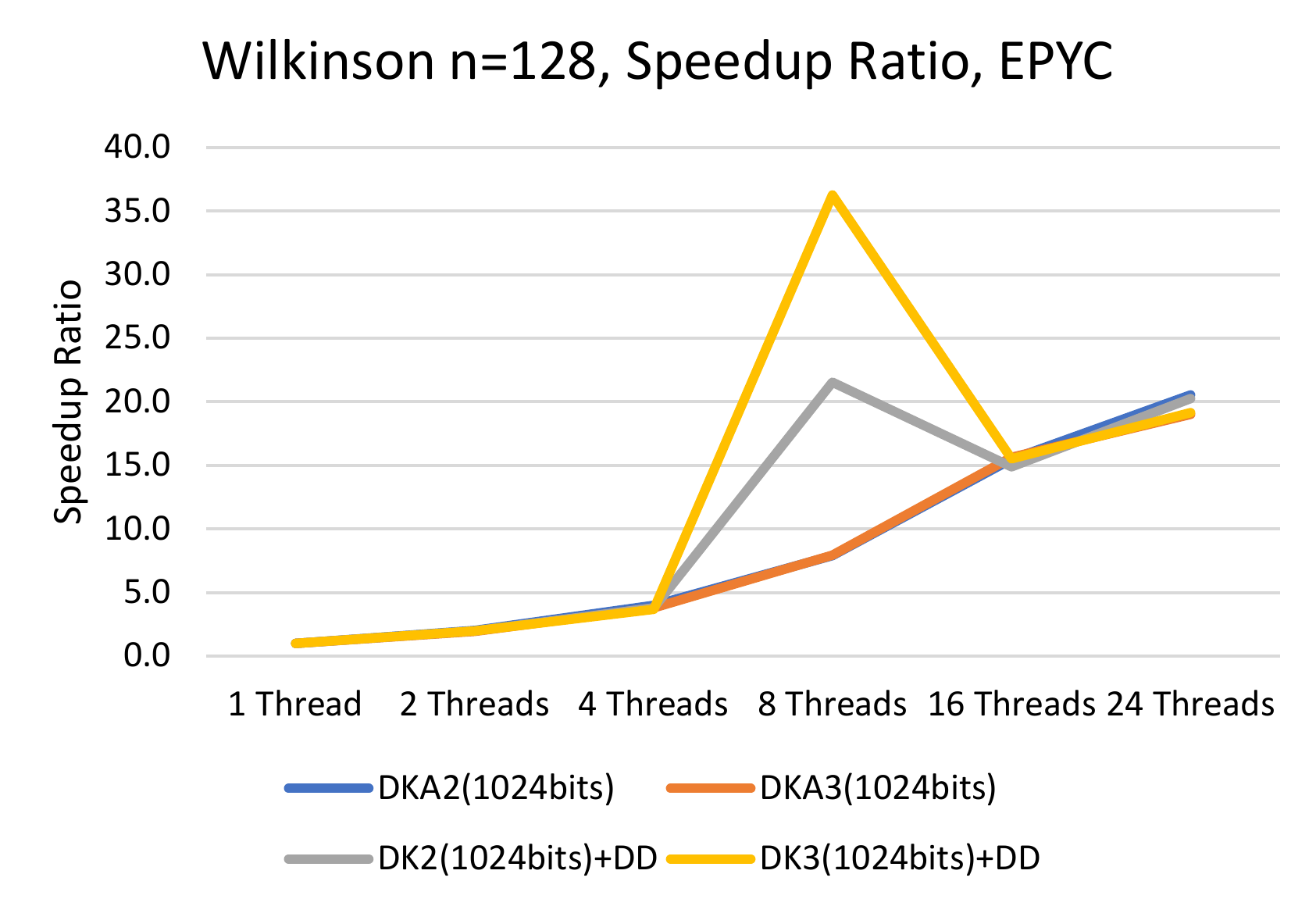}
		\includegraphics[width=.45\textwidth]{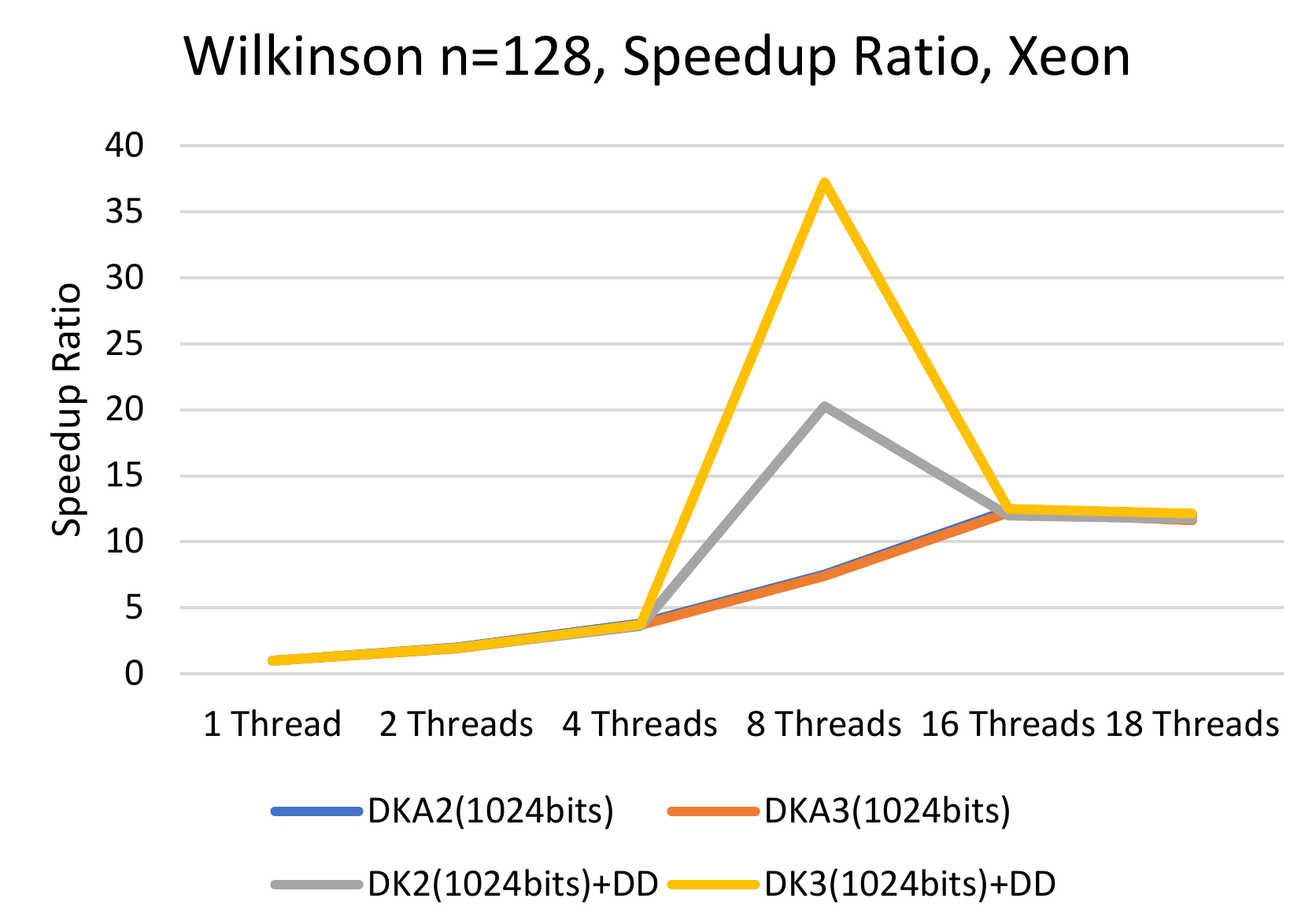}
		\caption{Speedup Ratio of Wilkinson's problem}\label{fig:speedup_ratio_wilkinson_128deg}
	\end{center}
\end{figure}

The speedup ratio by parallelization with OpenMP is illustrated in \figurename\ref{fig:speedup_ratio_wilkinson_128deg}. 

These figures illustrate that parallelization is efficient for direct iterative methods, and that the decrease of iterative times pull up the speedup ratio at DK2+DD and DK3+DD using eight threads.

%--------------------------------------%
% 
%--------------------------------------%
\subsection{Chebyshev quadrature problem}

Here, we explain the results of solving the Chebyshev quadrature formula quantile problem ($n = 256$ and $n = 512$).

\begin{figure}[htbp]
	\begin{center}
		\includegraphics[width=.55\textwidth]{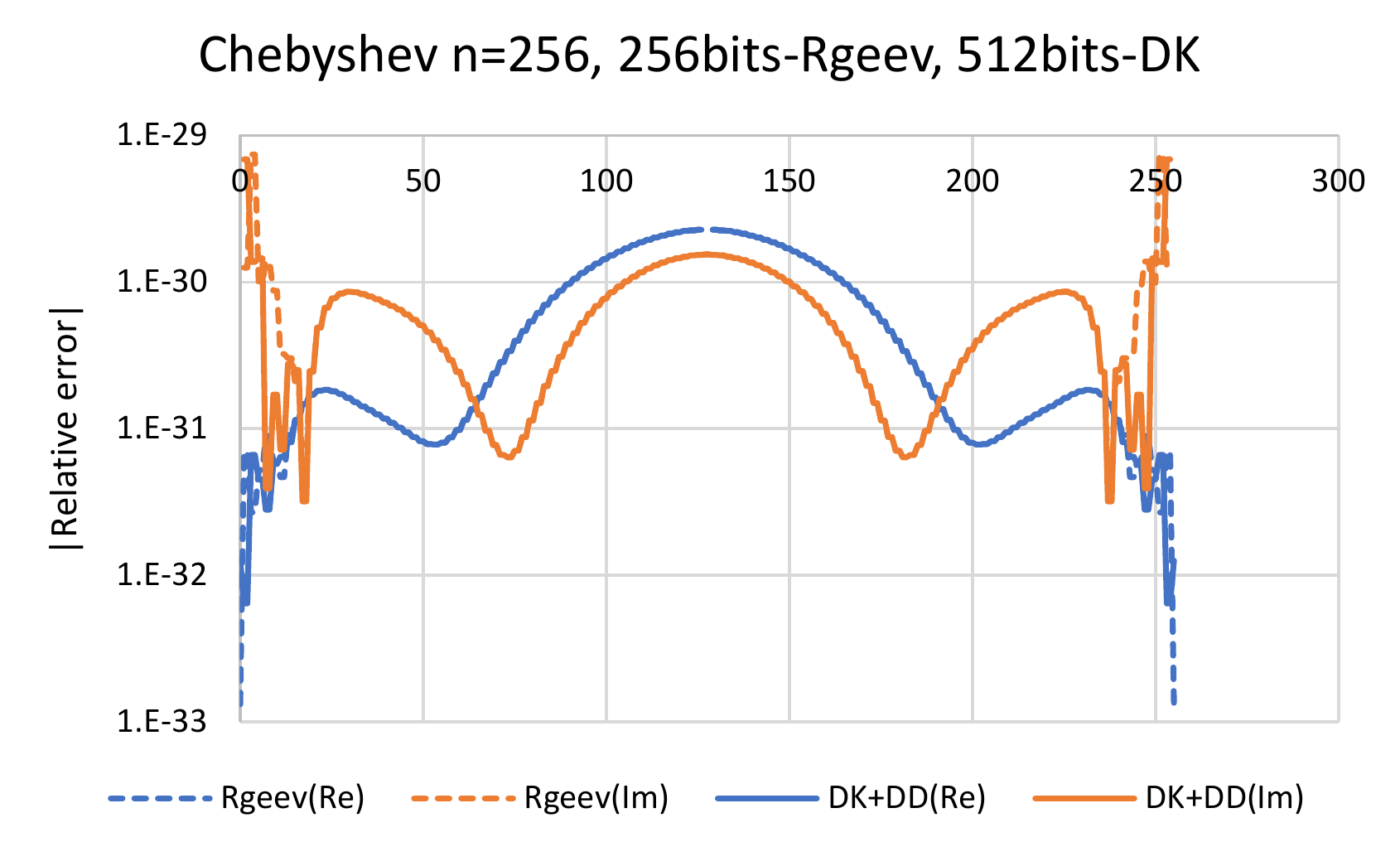}
		\caption{Relative errors of Rgeev(256 bits) and DK-DD(512 bits) for the Chebyshev quadrature problem: $n=256$.}\label{fig:relerr_cheb_256deg_rgeev256_dka512}
	\end{center}
\end{figure}

First we present the case of $n = 256$. MPREAL (256 bits) is adopted to derive the coefficients based on the (\ref{eqn:chebcoef}) formula, and the DK methods employs a 512 bits calculation and $\varepsilon_{\rm rel} := 8.6\times 10^{-68}$ and $\varepsilon_{\rm abs} := 1.0\times 10^{-300}$ for convergence determination. As illustrated in \figurename\ref{fig:relerr_cheb_256deg_rgeev256_dka512}, this results in accurate approximations of roots of approximately 30 decimal digits.

\tablename\ \ref{table:cheb256} presents all the results of the benchmark tests. MPREAL (256 bits) Rgeev took 57.6 s on EPYC and 38.4 on Xeon, respectively. The DK methods employ a 512 bits MPF arithmetic, $\varepsilon_{\rm rel} := 8.6\times 10^{-68}$ and $\varepsilon_{\rm abs} := 1.0\times 10^{-300}$ for convergence determination.

\begin{table}[htpb]
	\begin{center}\small
		\caption{Chebyshev quadrature :$n=256$, Rgeev 256 bits, DK methods 512 bits}\label{table:cheb256}
	\scalebox{0.8}{\begin{tabular}{|c|c|c|c|c|c|c|c|c|c|c|}
	\multicolumn{1}{c|}{EPYC} & \multicolumn{2}{|c|}{Rgeev (Second)} &\multicolumn{2}{|c|}{DKA2(512 bits)}	&\multicolumn{2}{|c|}{DKA3(512 bits)}	&\multicolumn{2}{|c|}{DK2(512 bits)+DD}	&\multicolumn{2}{|c|}{DK3(512 bits)+DD}	\\ \hline
	\#Thr.&	DD	&MPFR 256 bits 	& Second	&\#Iter.&Second	&\#Iter.&Second	&\#Iter.&Second	&\#Iter. \\ \hline
	1 &	1.29	&57.6			&203.0		&1344	&183	&671	&\underline{47}		&316	&70.3	&264 \\
	2 &			&				&100.0		&1344	&93.8	&671	&\underline{7.11}	&96		&\underline{10}		&73	\\
	4 &			&				&\underline{51.1}		&1344	&\underline{46.3}	&671	&\underline{3.62}	&96		&\underline{5.08}	&71	\\
	8 &			&				&\underline{25.3}		&1344	&\underline{23.1}	&671	&\underline{5.77}	&308	&\underline{8.92}	&264 \\
	16&			&				&\underline{13}			&1344	&\underline{11.7}	&671	&\underline{0.918}	&95		&\underline{1.22}	&71 \\
	24&			&				&\underline{8.96}		&1344	&\underline{8.1}	&671	&\underline{2.02}	&263	&\underline{3.1}	&263 \\ \hline
	\end{tabular}}	
	\scalebox{0.8}{\begin{tabular}{|c|c|c|c|c|c|c|c|c|c|c|}
	\multicolumn{1}{c|}{Xeon} & \multicolumn{2}{|c|}{Rgeev (Second)} &\multicolumn{2}{|c|}{DKA2(512 bits)}	&\multicolumn{2}{|c|}{DKA3(512 bits)}	&\multicolumn{2}{|c|}{DK2(512 bits)+DD}	&\multicolumn{2}{|c|}{DK3(512 bits)+DD}	\\ \hline
	\#Thr.&	DD	&MPFR 256 bits 	& Second	&\#Iter.&Second	&\#Iter.&Second	&\#Iter.&Second	&\#Iter. \\ \hline
	1 &	1.12	&38.4		&146.0	&1344	&127	&671	&\underline{32.1}	&316	&49.2	&264 \\
	2 &			&			&69.0	&1344	&64.8	&671	&\underline{4.95}	&96		&\underline{6.91}	&73 \\
	4 &			&			&\underline{36.3}	&1344	&\underline{34.2}	&671	&\underline{2.61}	&96		&\underline{3.59}	&71 \\ 
	8 &			&			&\underline{18.6}	&1344	&\underline{17.7}	&671	&\underline{4.2}	&308	&\underline{6.63}	&264 \\
	16&			&			&\underline{10.9}	&1344	&\underline{10}		&671	&\underline{0.779}	&95		&\underline{1.07}	&71 \\
	18&			&			&\underline{10.8}	&1344	&\underline{9.79}	&671	&\underline{2.32}	&293	&\underline{3.62}	&250 \\ \hline
	\end{tabular}}	
	\end{center}
\end{table}

The common numerical properties and trend of computational times is presented in \tablename\ \ref{table:cheb256}. In addition, we observe that the serial DD2+DD is faster than MPREAL Rgeev.

\begin{figure}[htbp]
	\begin{center}
		\includegraphics[width=.45\textwidth]{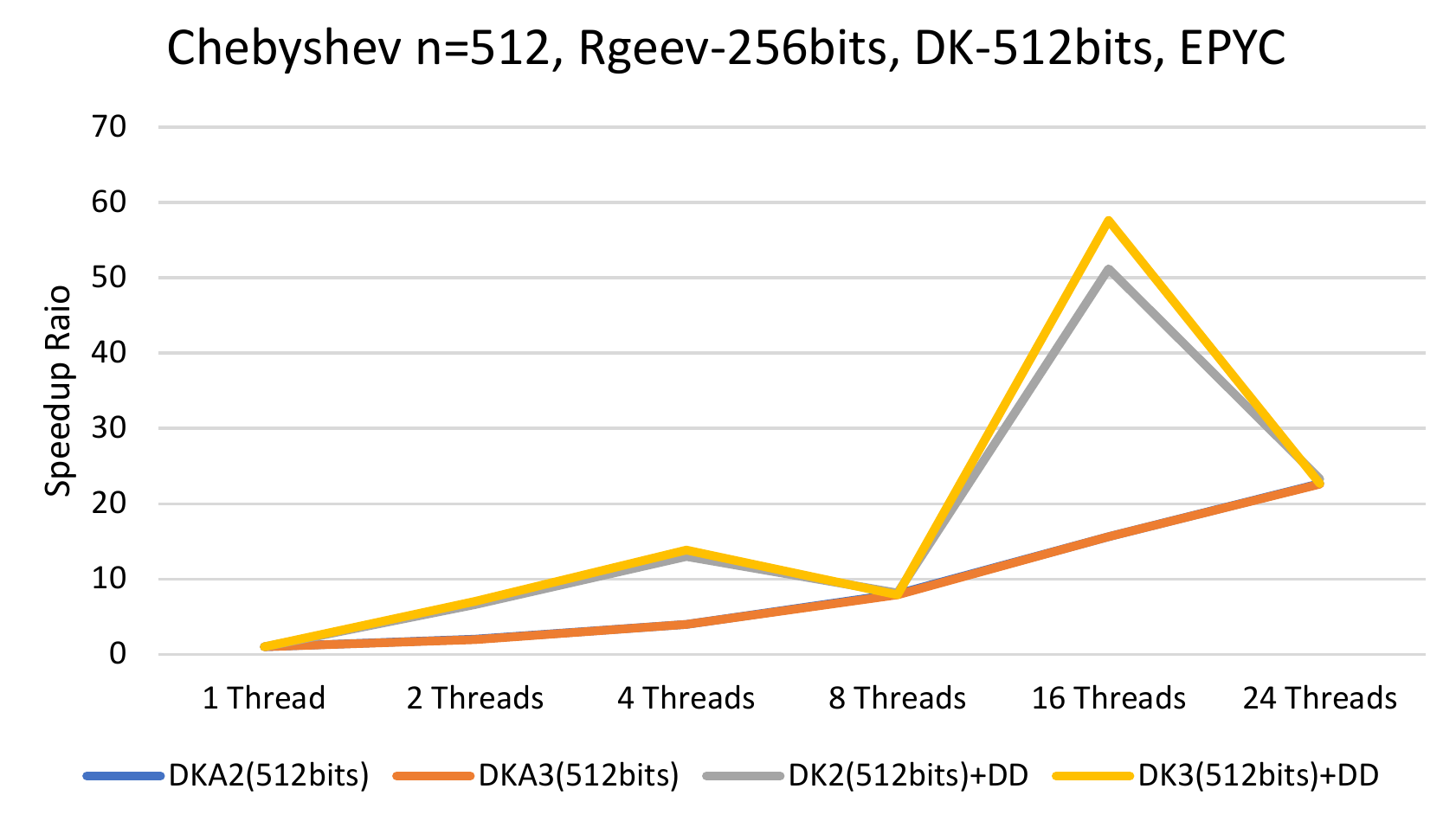}
		\includegraphics[width=.45\textwidth]{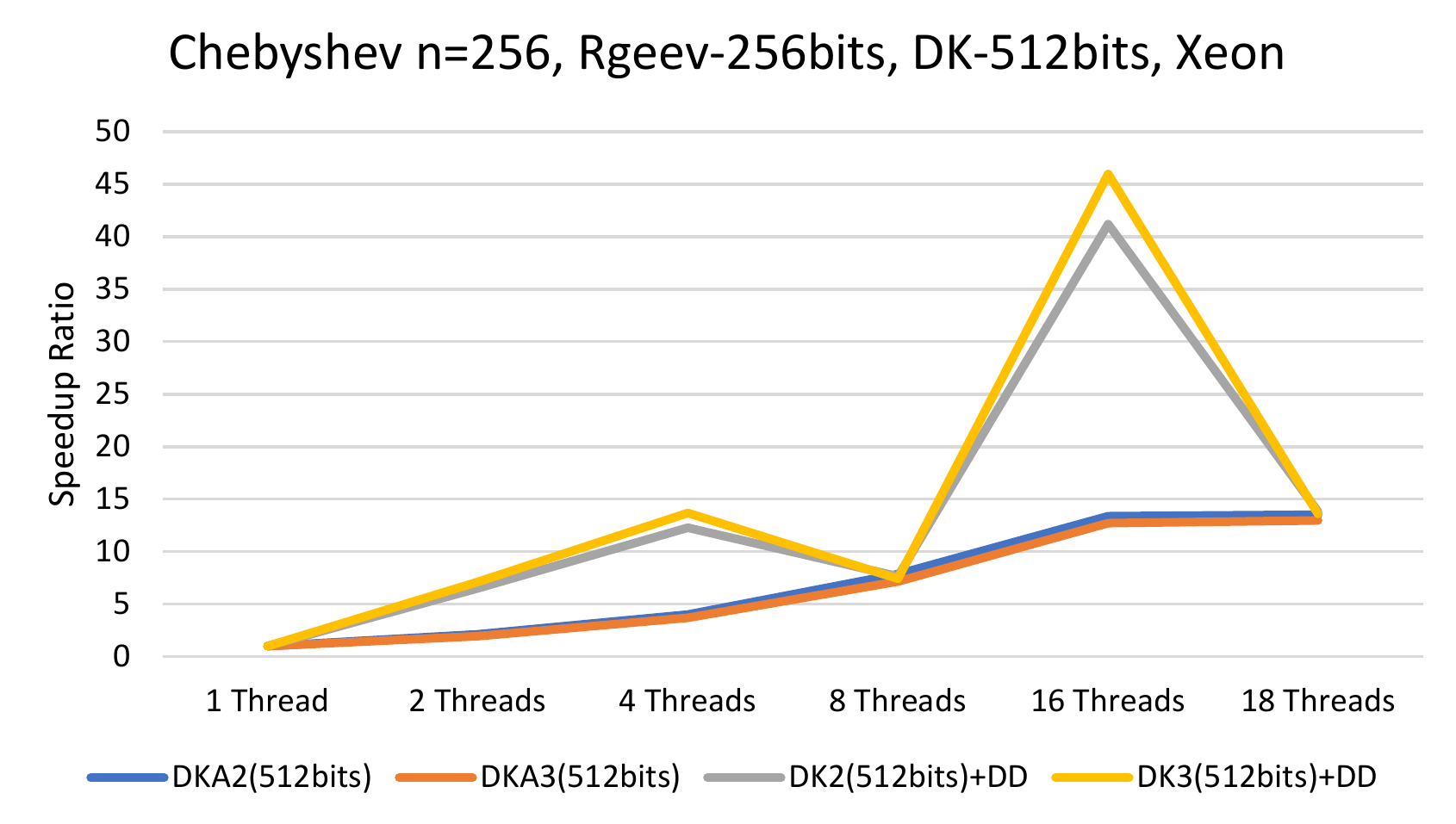}
		\caption{Speedup Ratio of Chebyshev quadrature: $n=256$}\label{fig:speedup_ratio_cheb_256deg}
	\end{center}
\end{figure}
	
According to the results presented in \tablename\ \ref{table:cheb256}, \figurename\ref{fig:speedup_ratio_cheb_256deg} illustrates the speedup ratio. We can confirm that second and third order DK methods have approximately achieved the ideal speedup ratio.

\begin{figure}[htbp]
	\begin{center}
		\includegraphics[width=.55\textwidth]{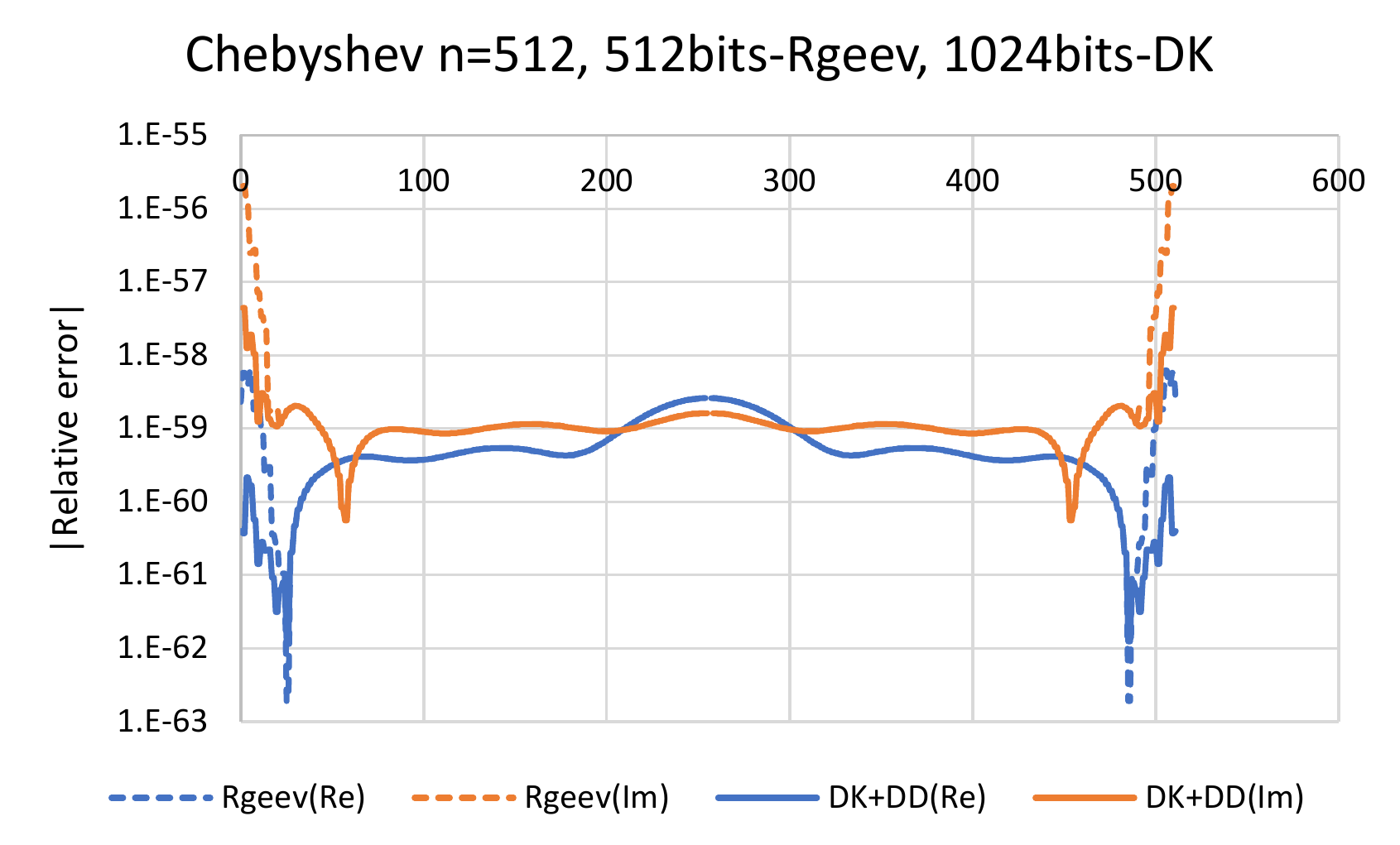}
		\caption{Relative errors of Rgeev(512 bits) and DK-DD(1024 bits) for the Chebyshev quadrature problem: $n=512$.}\label{fig:relerr_cheb_512deg_rgeev512_dka1024}
	\end{center}
\end{figure}

Second, the results of solving the Chebyshev quadrature problem in the case of $n = 512$ are indicated. MPREAL (512 bits) was adopted to derive the coefficients based on the (\ref{eqn:chebcoef}) formula, and the MPREAL 512 bits Rgeev took 636.0 s and 441.0 s on EPYC and Xeon, respectively.
The DK method employs a 1024 bits calculation, $\varepsilon_{\rm rel} := 7.5\times 10^{-145}$, and $\varepsilon_{\rm abs} := 1.0\times 10^{-300}$ for the convergence decision. This results in correct approximate solutions from 58 to 62 decimal digits, as illustrated in \figurename\ref{fig:relerr_cheb_512deg_rgeev512_dka1024}.

All results obtained through the benchmark test are presented in \tablename\ \ref{table:cheb512}.

\begin{table}[htpb]
	\begin{center}\small
	\caption{Chebyshev quadrature :$n=512$, Rgeev 512 bits, DK methods 1024 bits}\label{table:cheb512}
	\scalebox{0.8}{\begin{tabular}{|c|c|c|c|c|c|c|c|c|c|c|}
	\multicolumn{1}{c|}{EPYC} & \multicolumn{2}{|c|}{Rgeev (Second)} &\multicolumn{2}{|c|}{DKA2(1024 bits)}	&\multicolumn{2}{|c|}{DKA3(1024 bits)}	&\multicolumn{2}{|c|}{DK2(1024 bits)+DD}	&\multicolumn{2}{|c|}{DK3(1024 bits)+DD}	\\ \hline
	\#Thr.&	DD	&MPFR 512 bits 	& Second	&\#Iter.&Second	&\#Iter.&Second	&\#Iter.&Second	&\#Iter. \\ \hline
		1 	&9.84	&636.0	&2850.0	&3032	&2460	&1510	&\underline{500}	&552	&\underline{577}	&357 \\
		2 	&		&		&1370.0	&3032	&1210	&1510	&\underline{90.9}	&200	&\underline{67.8}	&82 \\
		4 	&		&		&689	&3032	&\underline{614}	&1510	&\underline{46.1}	&202	&\underline{34.5}	&83 \\
		8 	&		&		&\underline{344}	&3032	&\underline{307}	&1510	&\underline{62.8}	&550	&\underline{77.3}	&363 \\
		16	&		&		&\underline{173}	&3032	&\underline{174}	&1510	&\underline{11.8}	&205	&\underline{8.59}	&84 \\ 
		24	&		&		&\underline{120}	&3032	&\underline{106}	&1510	&\underline{8.47}	&213	&\underline{6.57}	&83 \\ \hline
	\end{tabular}}		

	\scalebox{0.8}{\begin{tabular}{|c|c|c|c|c|c|c|c|c|c|c|}
	\multicolumn{1}{c|}{Xeon} & \multicolumn{2}{|c|}{Rgeev (Second)} &\multicolumn{2}{|c|}{DKA2(1024 bits)}	&\multicolumn{2}{|c|}{DKA3(1024 bits)}	&\multicolumn{2}{|c|}{DK2(1024 bits)+DD}	&\multicolumn{2}{|c|}{DK3(1024 bits)+DD}	\\ \hline
	\#Thr.&	DD	&MPFR 512 bits 	& Second	&\#Iter.&Second	&\#Iter.&Second	&\#Iter.&Second	&\#Iter. \\ \hline
	1 &	9.25	&441.0	&2030.0	&3032	&1740	&1510	&\underline{362}	&552	&\underline{407}	&357 \\
	2 &			&		&1020.0	&3032	&937	&1510	&\underline{67.3}	&200	&\underline{48.2}	&82 \\
	4 &			&		&528	&3032	&471	&1510	&\underline{35.3}	&202	&\underline{25.5}	&83 \\
	8 &			&		&\underline{267}	&3032	&\underline{238}	&1510	&\underline{48.5}	&550	&\underline{56.5}	&363 \\
	16&			&		&\underline{159}	&3032	&\underline{142}	&1510	&\underline{10.8}	&205	&\underline{7.85}	&84 \\
	18&			&		&\underline{151}	&3032	&\underline{132}	&1510	&\underline{26.9}	&539	&\underline{31.7}	&365 \\ \hline
	\end{tabular}}		
	\end{center}
\end{table}

We confirm that more cases of DD2+DD and DD3+DD are faster than MPFR Rgeev. 

\begin{figure}[htbp]
	\begin{center}
		\includegraphics[width=.45\textwidth]{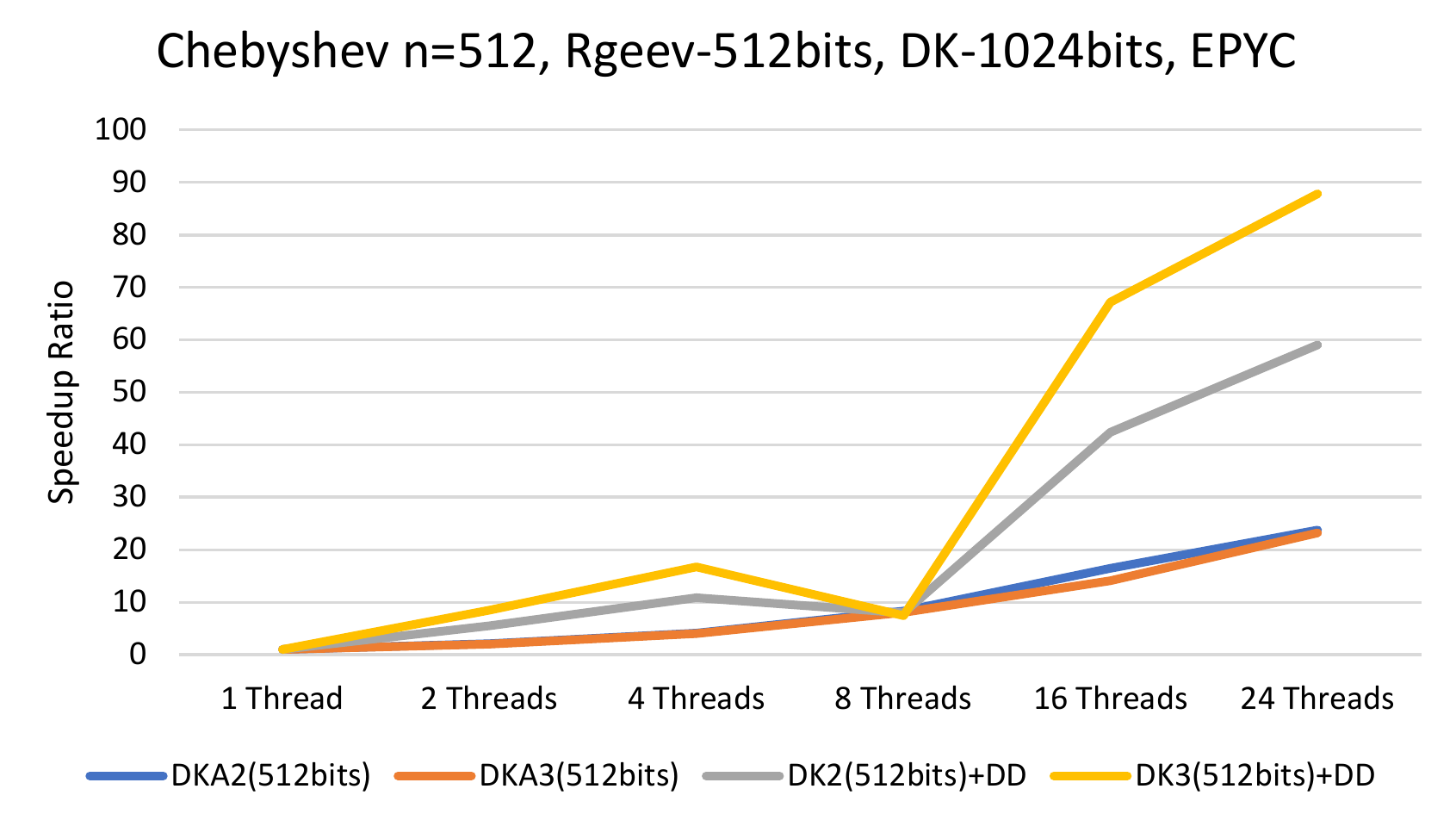}
		\includegraphics[width=.45\textwidth]{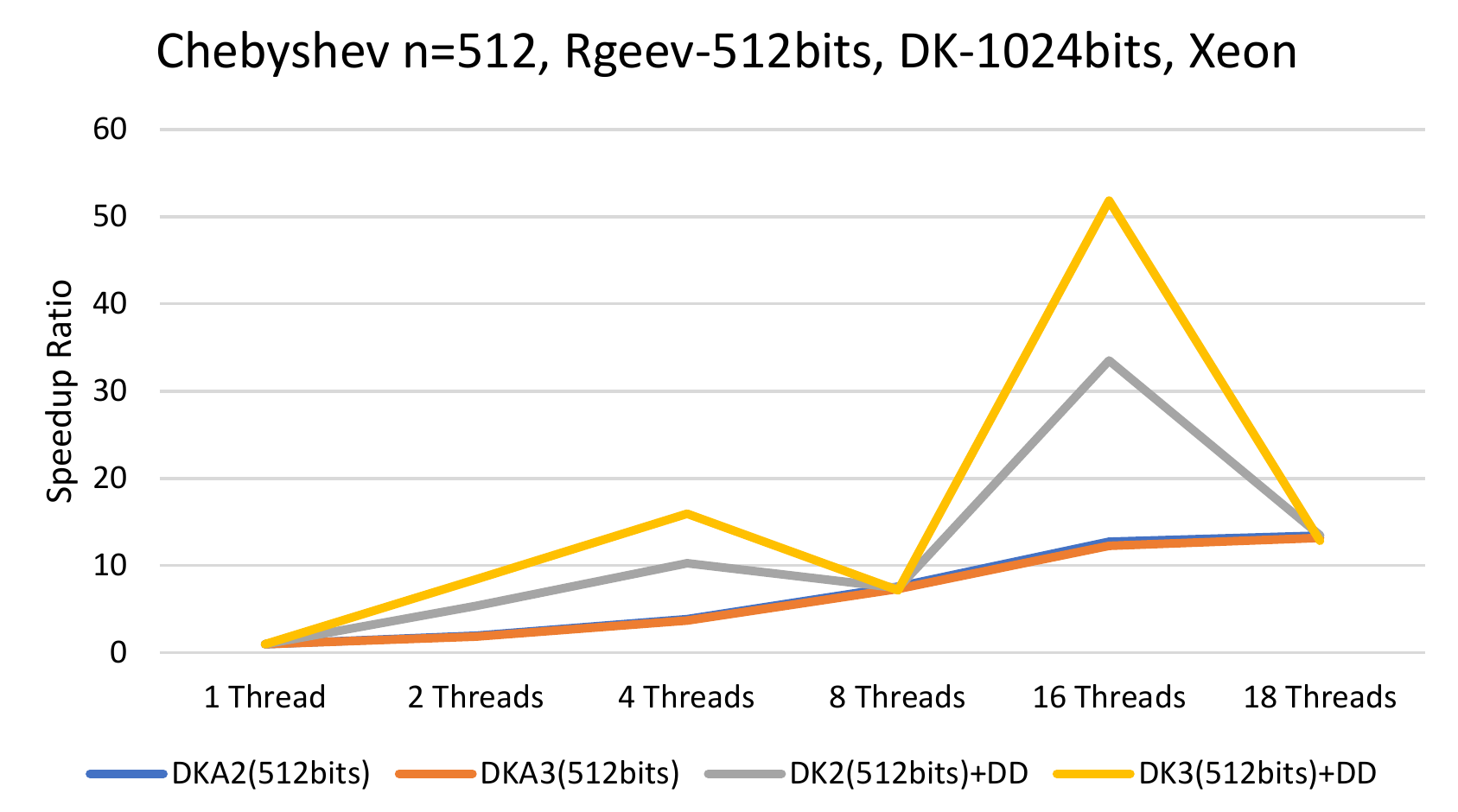}
		\caption{Speedup Ratio of Chebyshev quadrature: $n=512$}\label{fig:speedup_ratio_cheb_512deg}
	\end{center}
\end{figure}

The speedup ratio of parallelized DK methods is illustrated in \figurename\ref{fig:speedup_ratio_cheb_512deg}.

These two examples indicate that even if the low precision eigensolver cannot obtain good approximations, there are examples where speedup can be easily achieved using low-precision results as an initial guess for direct iterative methods with higher precision. Of course, it is not useful for all ill-conditioned algebraic equations, but the easy-to-implement and highly parallelizable direct iterative method may have a reasonably wide range of applications.

%--------------------------------------%
% 
%--------------------------------------%
\section{Conclusion and future works}

For two types of ill-conditioned algebraic equations, initial guesses were obtained using a low-precision eigensolver, and it was indicated that a highly accurate and efficient direct iterative method can be accelerated by parallel computing. Considering the stability of the iterative computing process, it is possible to speed up the eigensolver, for example, by using the approximate eigenvalues from low-precision computation as the origin shift in high-precision computation, but the direct iterative method for algebraic equations with its computational complexity, parallelism, and ease of implementation, is also useful in long precision computing environments. In several cases, it can be stated that our mixed precision approach is useful in solving high degree algebraic equations with long precision floating-point arithmetic. 

For future studies, we will implement and confirm the effectiveness of our approach as follows :
\begin{enumerate}
	\item Application to higher degree ill-conditioned algebraic equations,
	\item Implementation of third order direct iteration methods and comparison among them,
	\item Application to general linear equations with Chebyshev Proxy Rootfinder.
\end{enumerate}

In addition, we accelerate the above applications using multi-component MPF arithmetic.

%---------------------------------------%
% 
%---------------------------------------%
\section*{Acknowledgment}
This study was supported by JSPS KAKENHI, Grant Number JP20K11843, and Shizuoka Institute of Science and Technology. We acknowledge all organizations that are continuously encouraging our study.

%---------------------------------------%
% 
%---------------------------------------%
%\section*{References}
%\bibliographystyle{splncs04}
%\bibliography{tkouya_utf8}

%---------------------------------------%
% End of document
%---------------------------------------%
\end{document}